\documentclass{article}

\usepackage{arxiv}

\usepackage{setspace}
\usepackage{lmodern}

\usepackage{textgreek}

\usepackage{fixltx2e}

\reversemarginpar

\usepackage{amsmath,amssymb,amsfonts,amsthm}
\usepackage{mathrsfs}

\usepackage[locale = FR]{siunitx} 
\DeclareSIUnit\vitesse{\meter\per\second}
\usepackage{eurosym}
\DeclareSIUnit{\octet}{o}

\usepackage{graphicx,array,tikz,multirow}
\usepackage{caption,subcaption} 
\usepackage{svg,float}
\usepackage{booktabs,paralist}
\newcolumntype{x}[1]{>{\centering\arraybackslash\hspace{0pt}}p{#1}}
\usepackage[section]{placeins}
\usepackage{hanging}

\usepackage{fancyhdr,emptypage} 

\fancypagestyle{plain}{ 
    \fancyhead{}\fancyfoot[C]{\thepage}

}

\usepackage[Conny]{fncychap}

\usepackage[francais,nohints,tight]{minitoc}		

\usepackage[nottoc]{tocbibind} 
\usepackage[titles]{tocloft}

\usepackage{textcomp} 

\usepackage{titling}
\usepackage{lipsum} 
\usepackage{csquotes} 

\usepackage{xspace}
\usepackage{afterpage}

\usepackage[textsize=footnotesize]{todonotes}

\usepackage{bookmark}
\usepackage{acronym}
\usepackage[nameinlink,french]{cleveref} 
\Crefname{figure}{Fig.}{Figs.} 
\crefname{figure}{fig.}{figs.}
\Crefname{equation}{Eq.}{Eqs.}
\crefname{equation}{eq.}{eqs.}
\Crefname{table}{Table.}{Tables.}
\crefname{table}{table.}{tables.}

\definecolor{color_ref}{rgb}{1.0, 0.13, 0.32} 
\definecolor{color_link}{rgb}{0.0, 0.0, 1.0}
\definecolor{curcolor}{rgb}{0.0, 0.0, 1.0} 
\definecolor{brightpink}{rgb}{1.0, 0.0, 0.5} 
\definecolor{navyblue}{rgb}{0.0, 0.0, 1.0}
\hypersetup{
	colorlinks=true, 
	pdfstartview=FitV, 
	urlcolor=brightpink, 
	linkcolor= navyblue, 
	citecolor=color_ref 
}

\usepackage[utf8]{inputenc} 
\usepackage[T1]{fontenc}    
\usepackage{hyperref}       
\usepackage{url}            
\usepackage{booktabs}       
\usepackage{amsfonts}       
\usepackage{nicefrac}       
\usepackage{microtype}      
\usepackage{lipsum}
\usepackage{graphicx}
\usepackage{xcolor}
\def\grad{\mbox{\rm{grad}}} 
\graphicspath{ {./images/} }
\usepackage{amsmath}
\theoremstyle{definition}
\newtheorem{theorem}{Theorem}[section]

\newtheorem{definition}[theorem]{Definition}

\newtheorem{exemple}[theorem]{Example}

\newtheorem{propriete}[theorem]{Propriety}

\newtheorem{lemma}[theorem]{Lemma}

\newtheorem{notation}[theorem]{Notations}
\newtheorem{problem}[theorem]{Problem}

\newtheorem{proposition}[theorem]{Proposition}
\newtheorem{remark}[theorem]{Remark}

\title{A new geometric approach for sensitivity analysis in linear programming}
\author{
 Mustapha Kaci$^{{\color{white}.}\bf 1}$ and Sonia Radjef$^{{\color{white}.}\bf 2}$\\
 $^{\bf 1,2{\color{white}.}}$Department of Mathematics, Signal Image Parole (SIMPA) Laboratory\\ University of Oran Mohamed Boudiaf USTO-MB, Oran, Algeria.\\
  {\color{blue}\texttt{kaci.mustapha.95@gmail.com} }$^{{\color{white}.}\bf 1}$\\
	{\color{blue}\texttt{soniaradjef@yahoo.fr} }$^{{\color{white}.}\bf 2}$
   }

\begin{document}
\fontfamily{ptm}\selectfont
\maketitle
\begin{abstract}
 In this paper, we present a new geometric approach for sensitivity analysis in linear programming that is computationally practical for a decision-maker to study the behavior of the optimal solution of the linear programming problem under changes in program data. First, we fix the feasible domain (fix the linear constraints). Then, we geometrically formulate a linear programming problem. Next,  we give a new equivalent geometric formulation of the sensitivity analysis problem using notions of affine geometry. We write the coefficient vector of the objective function in polar coordinates and we determine all the angles for which the solution remains unchanged. Finally, the approach is presented in detail and illustrated with a numerical example.
\end{abstract}

\section{Introduction}
 \hspace{0.38cm}\par One of the star models of operational research is linear programming. It is a field of mathematical programming the most studied. It concerns the optimization of a mathematical program where the objective function and the functions defining the constraints are linear, see \cite{ref2,ref4}. Geometrically, linear programming problems are of convex programming, the linear constraints form a convex polyhedron; so, the results of the convexity are exploited. The vertices of the convex polyhedron form the basic solutions, one of them can be the optimal solution.
\par One of the essential parts of linear programming is sensitivity analysis, also called post-optimality analysis, because it starts from the original optimal solution, see \cite{ref1,ref3}.\\
Linear programming  models concrete problems such as maximizing a company's profits, but changes in market data require updates for the initial problems, and the sensitivity analysis is used to illustrate the margins of the linear program parameters for which the solution of the initial problem remains stable. Another situation that requires a stability study of the optimal solution is the committed errors in the data of the problem and it will allow us to avoid  restarting the procedure of resolution several times.
\par
In this work, we give a geometric approach to study the stability of the optimal solution when we change the coefficients of two decision variables in the objective function. This approach allows us to do the stability study with the variation of the two parameters  at the same time. It uses the Euclidean norm as well as the linearity of the constraint functions to formulate a new problem  equivalent to sensitivity analysis problem, then this new problem is solved with a simple calculation.   We fix a corner of the feasible domain and we find all the linear forms that reach their optimum value at this point. Geometrically, this consists to determine hyperplanes   in $ \mathbb{R}^{3}$, therefore the variation of the coefficients of the objective function corresponds to the rotation of the plan.\\
First, we consider a linear  maximization programming problem (\ref{Pbm0}), then we give a geometric formulation of a linear programming problem. Then, we use the latter to determine an equivalent problem to sensitivity analysis problem which consists in determining an angle interval on which the vector of the objective function coefficients is allowed to rotate to determine a linear form that reached its optimal value in initial optimal solution. Then, we look to the objective function  as a plan  of $\mathbb{R}^{3}$, and the linear constraints as a convex subset  of $\mathbb{R}^{2}$. The variation of the objective function coefficients corresponds to the rotation of the plan.
\par This paper is structured as follow: First, we  give some preliminaries in section \ref{section 2}. Then, we give the mathematical formulation of a sensitivity analysis problem in section \ref{section 3}.  Next, a  geometric formulation of a linear programming problem is given in section \ref{section 4}. After that, we give a geometric approach to sensitivity analysis in linear programming  in section \ref{section 5}. Then, we illustrate this approach with a numerical example in section \ref{section 6}. Finally, a conclusion is given in section \ref{section7}.

\section{Preliminary}\label{section 2}
\begin{figure}
\begin{center}
\includegraphics[scale=0.5]{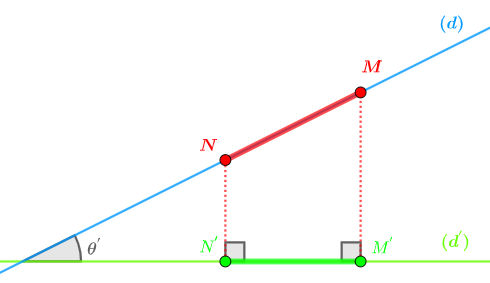}
\caption{\label{im36}Projection of a segment on a line.}

\end{center}
\end{figure}

\begin{notation}
Let  $A=(a_{ij})$ be a matrix of size $m_{1}\times m_{2}$, where $m_{1}, m_{2}$ are any integers. We  denote by $A^{t}$ the transposed matrix of $A$, defined by
$$
A^{t}=(a_{ji})\,.
$$
\end{notation}
\begin{definition}(Isometry of $\mathbb{R}^{2}$) We call an isometry of  $\mathbb{R}^{2}$, a linear application $u\in\mathcal{L}\left(\mathbb{R}^{2},\mathbb{R}^{2}\right)$ which checks one of the following properties:
\begin{itemize}\itemsep2pt \parskip2pt \parsep2pt
\item[a)] $u$ preserves the norm: $\forall x\in\mathbb{R}^{2}$,\hspace{0.15cm}$\left\|u(x)\right\|=\left\|x\right\|$,
\item[b)] $u$ preserves the scalar product: $\forall x,y\in\mathbb{R}^{2}$,\hspace{0.15cm}$\left\langle u(x),u(y)\right\rangle=\left\langle x,y\right\rangle$,
\item[c)] $u$ transforms an orthonormal basis into an orthonormal basis.
\end{itemize}
The application $u$ is also called orthogonal automorphism. The vector space of all the isometries (orthogonal automorphisms) of $\mathbb{R}^{2}$ is denoted by $\mathcal{O}\left(\mathbb{R}^{2}\right)$.
\end{definition}
\begin{propriete} Let $u_{1},u_{2}$ be two isometries of $\mathcal{O}\left(\mathbb{R}^{2}\right)$, then
\begin{enumerate}\itemsep1pt \parskip1pt \parsep1pt
\item $u_{1}$ and $u_{2}$ are one to one maps.
\item  $u_{1}\circ u_{2}$ and $u^{-1}_{1}$ are  isometries.
\item The translations and rotations of  $\mathbb{R}^{2}$ are  isometries.
\end{enumerate}
\begin{exemple} $\;$\\
\begin{enumerate}\itemsep1pt \parskip1pt \parsep1pt
\item Let $v$ be a vector in $\mathbb{R}^{2}$. A translation of $\mathbb{R}^{2}$ in the direction of $v$ is a one to one application noted by $T_{v}$, which associates to an element $x$ of $\mathbb{R}^{2}$ another element $T_{v}(x)=x+v$ of $\mathbb{R}^{2}$.
\item Let $\theta$ be an angle between $0$ and $2\pi$.  A rotation of $\mathbb{R}^{2}$ of angle $\theta$ is a one to one application, defined by
$$
\begin{array}{cccc}
R_{\theta}:&\mathbb{R}^{2}&\longrightarrow&\mathbb{R}^{2}\\
           &(x_{1},x_{2}) &\longrightarrow&R_{\theta}(x_{1},x_{2})=\text{Rot}(\theta)(x_{1},x_{2})^{t},
\end{array}
$$
where $\text{Rot}(\theta)$ is the matrix associated to $R_{\theta}$, defined by
$$
\text{Rot}(\theta)=\left(\begin{array}{cr}\cos(\theta)&-\sin(\theta)\\
                                \sin(\theta)&\cos(\theta)
               \end{array}\right).
$$
\end{enumerate}
Note that for any angle $\theta$, we have
\begin{align*}
\text{Rot}(\theta)^{t}\text{Rot}(\theta)&=\left(\begin{array}{rc}\cos(\theta)&\sin(\theta)\\
                                -\sin(\theta)&\cos(\theta)
                              \end{array}\right)\left(\begin{array}{cr}\cos(\theta)&-\sin(\theta)\\
                                \sin(\theta)&\cos(\theta)
                              \end{array}\right)\\
            &=\left(\begin{array}{ccc}1&0\\
                                0&1
                    \end{array}\right)=\text{Id}_{2}.
\end{align*}
\end{exemple}
\begin{definition}
 Let $F_{1}$ and $F_{2}$ be two supplementary linear subspaces of $\mathbb{R}^{n}$. Then, we call a projection on $F_{1}$ parallel to $F_{2}$, the map $P_{F_{1}}$ which associates to $x$ of $\mathbb{R}^{n}$ the unique element $P_{F_{1}}(x)$ of $F_{1}$ such that $x=P_{F_{1}}(x)+z$ where $z\in F_{2}$. Such an application is also called a projector.
\end{definition}
\end{propriete}
\begin{propriete}\label{prorto}
Consider two lines ($d$) and ($\acute{d}$) that form an angle $\acute{\theta}$. If $M$ and $N$ are two points that belong to $(d)$ and $\acute{M}=P_{(d)}(M)$, $\acute{N}=P_{(\acute{d})}(N)$, their  orthogonal projection respectively, see the figure (\ref{im36}), we obtain
$$
\acute{M}\acute{N}= MN\cos(\acute{\theta}).
$$
\end{propriete}

\section{Problem formulation }\label{section 3}

Consider the following initial linear programming problem in the standard form:
\begin{equation}\label{Pbm0}
\left\{
\begin{array}{l}
      \underset{x}\max\hspace{0.15cm} f(x)=c_{0}^{t}x           \\
      Ax\leq b\\
x\geq0
\end{array},
\right.
\end{equation}
with
$$
\begin{array}{cccc}
A=\left(\begin{array}{ccccccccc}
a_{11}&a_{12}\\
a_{21}&a_{22}\\
 \vdots    & \vdots \\
a_{m1}&a_{m2}
\end{array}\right),&b=\left(\begin{array}{c}b_{1}\\b_{2}\\\vdots\\b_{m}\end{array}\right),&c_{0}=\left(\begin{array}{c}c^{0}_{1}\vspace{0.15cm}\\c^{0}_{2}\end{array}\right),&x=\left(\begin{array}{c}x_{1}\\x_{2}\end{array}\right).
\end{array}
$$
Where  $c_{0}$ is the vector of the objective function,  $c^{0}_{1},c^{0}_{2}$ are constants, $x$ is a $2\times1$ vector of decision variables, $A$ is a $m\times 2$ matrix of constants,  $b$ is a $m\times 1$ vector of constants  and $m$ is the number of linear constraints.
\par
With the hypothesis that $c_{1}^{0}, c_{2}^{0}$ are positive numbers, we will see in the following that it is an artificial hypothesis that we made only to simplify the continuation of this presentation, see the remark (\ref{reak}).
\begin{notation}
Let us define:
\begin{itemize}\itemsep2pt \parskip2pt \parsep2pt
\item The feasible region:
$$
S:=\left\{x\in\mathbb{R}^{2}\hspace{0.1cm} :\hspace{0.1cm} Ax\leq b,\hspace{0.15cm}  x\geq0 \right\}.
$$
\item The vector space $\mathcal{L}$  of linear formes defined on $\mathbb{R}^{2}$, is defined by:
$$
\mathcal{L}\left(\mathbb{R}^{2}\right):=\left\{g:\mathbb{R}^{2}\rightarrow\mathbb{R}\hspace{0.1cm}:\hspace{0.1cm}g(x_{1},x_{2})=c_{1}x_{1}+c_{2}x_{2}\right\} .
$$
\item $\left\|.\right\|$ denotes the Euclidean norm on $\mathbb{R}^{2}$.
\end{itemize}
\end{notation}

Let $x^{0}=(x^{0}_{1},x^{0}_{2})^{t}$ be an optimal solution of the problem (\ref{Pbm0}). A sensitivity analysis problem is to find all linear forms $g$ different from $f$  verifying
\begin{equation}\label{Pbm1}
\begin{array}{cc}
\arg  \underset{x\in S} \max \, f(x)= \arg \underset{x\in S} \max \, g(x),&  g\in\mathcal{L}\left(\mathbb{R}^{2}\right),
\end{array}
\end{equation}
where
$$
\arg  \underset{x\in S} \max \, f(x)=\left\{y\in S\hspace{0.2cm}|\hspace{0.2cm}f(y)\geq f(x),\hspace{0.2cm} x\in S\right\}.
$$
\section{Geometric formulation of a linear programming problem}\label{section 4}
\par
In this part, we will determine a problem equivalent to the linear programming problem (\ref{Pbm0}). In $\mathbb{R}^{3}$, we find that the feasible region $S$ is a polygon entirely contained in the vector hyperplane of equation $z=0$, we can easily show that the graph of $f$ is a vector hyperplane (linear subspace of dimension two).\\

 Let's start by transforming the expression of the linear form $f$ as follows:
\[
\begin{array}{ccc}
 f(x_{1},x_{2})=c^{0}_{1}x_{1}+c^{0}_{2}x_{2}&\Leftrightarrow&  c^{0}_{1}x_{1}+c^{0}_{2}x_{2} -f(x_{1},x_{2})=0\vspace{0.2cm}\\
                             &\Leftrightarrow&\left\langle (c^{0}_{1},c^{0}_{2},-1),(x_{1},x_{2},f(x_{1},x_{2}))\right\rangle=0 \vspace{0.2cm}\\
														 &\Leftrightarrow&\left\langle n_{f},(x_{1},x_{2},z)\right\rangle=0,
\end{array}
\]
where $x_{1},x_{2}\in\mathbb{R}$, $z=f(x_{1},x_{2})$  and $n_{f}=(c^{0}_{1},c^{0}_{2},-1)$.
So, the graph of the linear form $f$, defined by
	$$
	\begin{array}{ccl}
	Gf(f)&=&\left\{(x_{1},x_{2},z)\in\mathbb{R}^{3}\hspace{0.15cm}|\hspace{0.15cm}z=f(x_{1},x_{2})\right\}\vspace{0.2cm}\\
	     &=&\left\{(x_{1},x_{2},z)\in\mathbb{R}^{3}\hspace{0.15cm}|\hspace{0.15cm}\left\langle n_{f},(x_{1},x_{2},z)\right\rangle=0\right\}
	\end{array}
	$$
is orthogonal to $\text{vect}\left(n_{f}\right)$, the linear subspace generated by the vector $n_{f}$, and therefore any vector of $\mathbb{R}^{3}$ can be written as the sum of a vector of $Gf(f)$ and another vector of $\text{vect}\left(n_{f}\right)$.

Now, we position ourselves on $\mathbb{R}^{2}$ (the plane $z=0$)  and we consider the line $d_{0}$ defined by the intersection of $Gf(f)$ with the plane $z=0$:
\begin{equation}
(d_{0}): \, \,  \, \, \,
 c^{0}_{1}x_{1}+c^{0}_{2}x_{2}=0,
\end{equation}
$d_{0}$ is a linear subspace of  $\mathbb{R}^{2}$ of dimension 1 directed by the vector:
$$
v_{0}=(-c^{0}_{2},c^{0}_{1}).
$$

\begin{definition}
The graph of the line $d_{0}$, defined by:
$$
D_{0}=\left\{(x_{1},x_{2})\in\mathbb{R}^{2}\hspace{0.15cm} |\hspace{0.15cm} c^{0}_{1}x_{1}+c^{0}_{2}x_{2}=0 \right\}.
$$
The upper half-space of $\mathbb{R}^{2}$ deprived of ($d_{0}$), defined by:
$$
D_{0}^{+}=\left\{(x_{1},x_{2})\in\mathbb{R}^{2} \hspace{0.15cm}|\hspace{0.15cm} c^{0}_{1}x_{1}+c^{0}_{2}x_{2}>0 \right\}.
$$
The lower half-space of $\mathbb{R}^{2}$ deprived of $d_{0}$, defined by:
$$
D_{0}^{-}=\left\{(x_{1},x_{2})\in\mathbb{R}^{2} \hspace{0.15cm}|\hspace{0.15cm} c^{0}_{1}x_{1}+c^{0}_{2}x_{2}<0 \right\}.
$$
\end{definition}

\begin{lemma}\label{lemme1}
For all $x^{1}, x^{2}\in D_{0}^{+}$ (respectively  $y^{1}, y^{2}\in D_{0}^{-}$),  we have
\begin{equation}
\frac{f(x^{1})}{\left\|x^{1}-P_{(d_{0})}(x^{1})\right\|}=\frac{f(x^{2})}{\left\|x^{2}-P_{(d_{0})}(x^{2})\right\|},
\end{equation}
and
\begin{equation}
\frac{f(y^{1})}{\left\|y^{1}-P_{(d_{0})}(y^{1})\right\|}=\frac{f(y^{2})}{\left\|y^{2}-P_{(d_{0})}(y^{2})\right\|}.
\end{equation}
\end{lemma}

\begin{figure}
\begin{center}
\includegraphics[scale=0.4]{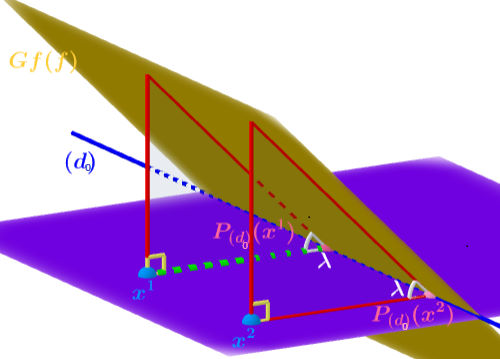}
\caption{\label{im33} Graphic illustration of the lemma  (\ref{lemme1}).}
\end{center}
\end{figure}

\begin{proof} The intersection of $Gf(f)$ with the plane $z=0$ forms a constant angle $0\leq\lambda\leq\pi$, see figure (\ref{im33}). Then,
$$
\left\{\begin{array}{lccccc}
\tan(\lambda)&=&\frac{f(x^{1})}{\left\|x^{1}-P_{(d_{0})}(x^{1})\right\|}&=&\frac{f(x^{2})}{\left\|x^{2}-P_{(d_{0})}(x^{2})\right\|},&x^{1}, x^{2}\in D_{0}^{+}\vspace{0.25cm}\\
\tan(\pi-\lambda)&=&\frac{f(y^{1})}{\left\|y^{1}-P_{(d_{0})}(y^{1})\right\|}&=&\frac{f(y^{2})}{\left\|y^{2}-P_{(d_{0})}(y^{2})\right\|},&y^{1}, y^{2}\in D_{0}^{-}
\end{array}.\right.
$$\end{proof}

\begin{remark}\label{reak}
If $c_{1}^{0}, c_{2}^{0}$ are any real numbers, we can come back to the hypotheses of the problem (\ref{Pbm0}) by making the following changes:

\begin{enumerate}\itemsep2pt \parskip2pt \parsep2pt
\item Apply a rotation $R_{\theta}$ over the feasible region $S$ and the vector $c_{0}$ simultaneously to obtain the problem:
\begin{equation}\label{Pbm0Rot}
\underset{x\in S}{\max}\hspace{0.15cm} \left(R_{\theta}(c_{0})\right)^{t}R_{\theta}(x),
\end{equation}
which is equivalent to the initial problem. Indeed, we have
\begin{align*}
\left(\text{Rot}(\theta)c_{0}\right)^{t}\text{Rot}(\theta)x&=c_{0}^{t}\left(\text{Rot}(\theta)^{t}\text{Rot}(\theta)\right)=c_{0}^{t}\,\mbox{\rm{Id}}_{2}x=c_{0}^{t}x=f(x).\\
\intertext{Then, an optimal solution $x^{0}$  of the problem (\ref{Pbm0}) is also optimal to the problem (\ref{Pbm0Rot}), indeed}
\left(\text{Rot}(\theta)c_{0}\right)^{t}\text{Rot}(\theta)x^{0}&=f(x^{0})\geq f(x)=\left(\text{Rot}(\theta)c_{0}\right)^{t}\text{Rot}(\theta)x.\\ \intertext{After that, we pose}
\big(\text{Rot}(\theta^{0})c_{0}\big)&=\left(\left|c_{1}^{0}\right|,\left|c_{2}^{0}\right|\right)\geq 0,\\
\intertext{where $\theta^{0}$ is a solution of the following linear system:}
\left\{
\begin{array}{cc}
      \cos(\theta)c_{1}^{0}-\sin(\theta)c_{2}^{0}&=\left|c_{1}^{0}\right|\\
      \sin(\theta)c_{1}^{0}+\cos(\theta)c_{2}^{0}&=\left|c_{2}^{0}\right|
       \end{array}.
\right.
\end{align*}

\item Apply a translation $T_{v_{0}}$ over  $S$ in the direction of $\text{Rot}\left(\theta^{0}\right).\left(c_{1}^{0},c_{2}^{0}\right)^{t}$ because the objective function is constant in the latter. Indeed, for all $y\in T_{v_{0}}(S)$, we have
$$
\begin{array}{ccc}
\begin{array}{rcl}
f\left(y\right)=f\left(x+v_{0}\right)&=&f\left(x\right)+f\left(v_{0}\right)\vspace{0.15cm}\\
               &\leq& f\left(x^{0}\right)+f\left(v_{0}\right)\vspace{0.15cm}\\
							 &\leq&f\left(x^{0}+v_{0}\right)\vspace{0.15cm}\\
							&\leq&f\left(y^{0}\right)
\end{array},&	x\in S.
\end{array}					
$$
where
$$
\begin{array}{ccc}
y^{0}=x^{0}+v_{0},\hspace{0.2cm}&\text{and}\hspace{0.2cm}&v_{0}=\alpha \text{Rot}\left(\theta^{0}\right)\left(c_{1}^{0},c_{2}^{0}\right)^{t},
\end{array}		
$$
 $\alpha$ is a strictly positive number sufficiently large.
\end{enumerate}
\end{remark}
\begin{proposition}\label{proposition1}
Let $c_{1}^{0}, c_{2}^{0} \geq0$, then the problem (\ref{Pbm0}) is equivalent to:
\begin{equation}\label{Pbm2}
\underset{x\in S}\max \, \,\left\|x-P_{(d_{0})}(x)\right\|.
\end{equation}

\end{proposition}

\begin{proof} $\;$\\
\begin{itemize}\itemsep1pt \parskip1pt \parsep1pt
\item[$\diamond$] \hspace{0.2cm} 	Let $x^{0}$ be an optimal solution of the problem (\ref{Pbm0}), then
\begin{equation}\label{eq1}
\begin{array}{cc}
f(x^{0})\geq f(x),&\hspace{0.2cm}\forall x\in S,
\end{array}
\end{equation}
from lemma (\ref{lemme1}), and inequality (\ref{eq1}), we get
$$
\begin{array}{cc}
\frac{f(x^{0})}{\left\|x^{0}-P_{(d_{0})}\left(x^{0}\right)\right\|}=\frac{f(x)}{\left\|x-P_{(d_{0})}(x)\right\|}\leq\frac{f\left(x^{0}\right)}{\left\|x-P_{(d_{0})}(x)\right\|},&\hspace{0.2cm}\forall x\in S.
\end{array}
$$
Then, for all $x\in S$, we get
\begin{equation}\label{supo}
\left\|x^{0}-P_{(d_{0})}\left(x^{0}\right)\right\|\geq\left\|x-P_{(d_{0})}(x)\right\|,
\end{equation}
which proves that $x^{0}$ is an optimal solution of the problem (\ref{Pbm2}).
\item[$\diamond$] Conversely, let $x^{0}$ be an optimal solution of the problem (\ref{Pbm2}), then $x^{0}$ verifies the inequality (\ref{supo}), then 
$$
\begin{array}{cc}
\frac{1}{\left\|x-P_{(d_{0})}(x)\right\|}\geq\frac{1}{\left\|x^{0}-P_{(d_{0})}\left(x^{0}\right)\right\|},&\hspace{0.2cm}\forall x\in S.
\end{array}
$$
We multiply the two sides of inequality  by $f(x)$, then we apply lemma (\ref{lemme1}) to the left side, and get
$$
\frac{f\left(x^{0}\right)}{\left\|x^{0}-P_{(d_{0})}\left(x^{0}\right)\right\|}=\frac{f(x)}{\left\|x-P_{(d_{0})}(x)\right\|}\geq\frac{f(x)}{\left\|x^{0}-P_{(d_{0})}\left(x^{0}\right)\right\|},
$$
which proves that $x^{0}$ is an optimal solution of the problem (\ref{Pbm0}),
$$
 f\left(x^{0}\right)\geq f(x),\, \,  \text{for all } \hspace{0.2cm} x\in S.
$$
\end{itemize}
\end{proof}
\section{A geometric approach to sensitivity analysis in linear programming}\label{section 5}
\par
\begin{figure}
\begin{center}
\includegraphics[scale=0.5]{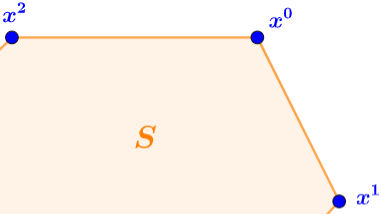}
\caption{\label{im31}Extreme points close to the optimal solution $x^{0}$.}
\end{center}
\end{figure}

 The geometric analysis we present here consists of writing a linear programming problem in the polar coordinate system and then determining a relationship that links the coefficients of the objective function $f$. We already know that to any linear form $g$ in $\mathcal{L}\left(\mathbb{R}^{2}\right)$, we can associate a single vector (gradient of $g$) in $\mathbb{R}^{2}$ and conversely. Another important fact is that we can get every corner of the feasible region  $S$ by solving a finite number of linear systems with two equations and two variables. This allows us to reformulate the problem (\ref{Pbm1}) as follows:
\begin{problem}\label{Pbm3}
Let $x^{1}, x^{0}, x^{2}$ be the successive corners (extreme points) of $S$, see figure (\ref{im31}), where $x^{0}$ is the optimal solution of the problem (\ref{Pbm2}). Find all  vector lines $(d)$ that simultaneously satisfy to the following  two  inequalities:
\begin{equation}\label{Pbm31}
                \left\|x^{0}-P_{(d)}\left(x^{0}\right)\right\|\geq\left\|x^{1}-P_{(d)}\left(x^{1}\right)\right\|,
\end{equation}
and
\begin{equation}\label{Pbm32}													
                \left\|x^{0}-P_{(d)}\left(x^{0}\right)\right\|\geq\left\|x^{2}-P_{(d)}\left(x^{2}\right)\right\| .	
   \end{equation}																
\end{problem}

\begin{proposition}
The problem (\ref{Pbm1}) and (\ref{Pbm3}) are equivalents.
\end{proposition}
\begin{proof}
This is a direct consequence of the proposition (\ref{proposition1}) and convexity of $S$.
\end{proof}

\begin{remark}\label{remarktra}
The fact that the problem (\ref{Pbm3}) is built from another problem that admits a solution, then it still admits the vector line $(d_{0})$ as a trivial solution.
\end{remark}

\begin{proposition}\label{solved}
Let $\theta_{1},\theta_{2}\in\left[0,\pi\right]$, $\phi\in\left[0,\frac{\pi}{2}\right[$, and $r,r_{1},r_{2}\geq0$, and
consider the following vectors $x^{10}$ and $x^{02}$ written in polar coordinate system:
$$
\begin{array}{crlcl}
x^{10}&:=&x^{1}-x^{0}&=&r_{1}(\cos(\theta_{1}),\sin(\theta_{1}))\vspace{0.2cm}\\
x^{02}&:=&x^{0}-x^{2}&=&r_{2}(\cos(\theta_{2}),\sin(\theta_{2}))\vspace{0.2cm}\\
c&=&(c_{1},c_{2})&=&r(\cos(\phi),\sin(\phi))\,.
\end{array}
$$
Then, the solutions of the problem (\ref{Pbm3}) are the line vectors defined by:
\begin{equation}
\begin{array}{cc}
(d) : r\cos(\phi)x_{1}+r\sin(\phi)x_{2}=0\quad with \hspace{0.25cm} \theta_{1}<\phi+\frac{\pi}{2}<\theta_{2}.
\end{array}
\end{equation}
\end{proposition}

\begin{figure}
\begin{center}
\includegraphics[scale=0.6]{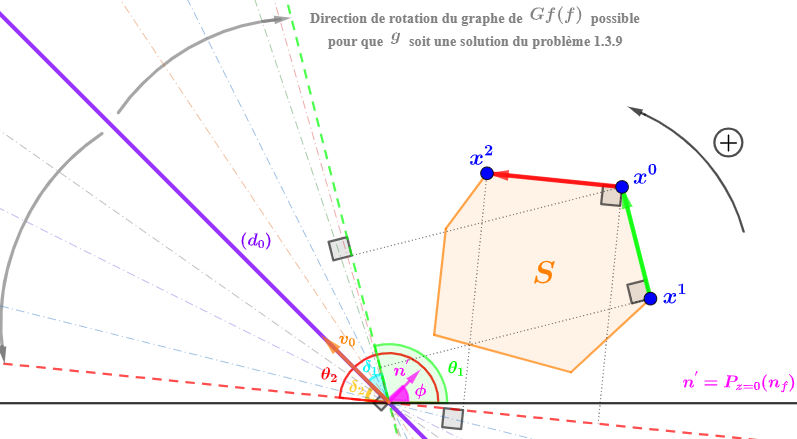}
\caption{\label{im32}Graphic illustration of the proposition (\ref{solved}).}
\end{center}
\end{figure}

\begin{proof} Let $x^{0}$ be a solution of the  problem (\ref{Pbm2}), then from propriety (\ref{prorto}), we get
$$
\left\{\begin{array}{l}
\left\|P_{(d)}(x^{1})-P_{(d)}(x^{0})\right\|=\left\|x^{1}-x^{0}\right\|\cos(\delta_{1})\vspace{0.25cm}\\
\left\|P_{(d)}(x^{2})-P_{(d)}(x^{0})\right\|=\left\|x^{2}-x^{0}\right\|\cos(\delta_{2})
\end{array},\right.
$$
where $\delta_{1}$ is the angle formed by the vector $x^{10}$ with $v_{0}$, and $\delta_{2}$ is the angle formed by the vector $x^{02}$ with $v_{0}$ (the order of vectors is important). On the other hand, see figure (\ref{im32}), we have
$$
\left\{\begin{array}{l}
                             |\delta_{1}|=|\phi+\frac{\pi}{2}-\theta_{1}|<\frac{\pi}{2}\vspace{0.2cm}\\
                             |\delta_{2}|=|\phi+\frac{\pi}{2}-\theta_{2}|<\frac{\pi}{2}
\end{array}.\right.
$$
Therefore, we get
$$
\left\{\begin{array}{l}
                          \left\|P_{(d)}(x^{1})-P_{(d)}(x^{0})\right\|<\left\|x^{1}-x^{0}\right\|\vspace{0.2cm}\\
                          \left\|P_{(d)}(x^{2})-P_{(d)}(x^{0})\right\|<\left\|x^{2}-x^{0}\right\|
\end{array}.\right.
$$
Since, for all $i=1, 2$ the quadrilateral $x^{0}x^{i}P_{(d)}(x^{i})P_{(d)}(x^{0})$ is irregular, the two segments $P_{(d)}(x^{0})x^{0}$ and $P_{(d)}(x^{i})x^{i}$ are parallels, then
$$
\left\{\begin{array}{l}
                          \left\|P_{(d)}(x^{1})-x^{1}\right\|<\left\|P_{(d)}(x^{0})-x^{0}\right\|\vspace{0.2cm}\\
                          \left\|P_{(d)}(x^{2})-x^{2}\right\|<\left\|P_{(d)}(x^{0})-x^{0}\right\|
\end{array}.\right.
$$

\end{proof}

\section{Numerical example}\label{section 6}
Consider the following linear programming problem:
\begin{equation} \label{s}
\left\{
\begin{array}{l}
      \underset{x_{1},x_{2}}\max\hspace{0.15cm}  f(x_{1},x_{2})=2x_{1}+3x_{2}           \\
      A(x_{1},x_{2})^{t}\leq b\\
x_{1}\geq 0; x_{2}\geq 0
\end{array},
\right.
\end{equation}
where
$$
\begin{array}{llll}
               A=\left(\begin{array}{ll}
                                  \frac{1}{4}&\frac{1}{2}\vspace{0.15cm}\\
					                        \frac{2}{5}&\frac{1}{5}\vspace{0.15cm}\\
																	0&\frac{4}{5}
																	
                           \end{array}\right),&\hspace{0.2cm}b= \left(\begin{array}{c}40\\40\\40\end{array}\right),&\hspace{0.2cm}c_{0}= \left(\begin{array}{c}2\\3\end{array}\right),&\hspace{0.2cm}x= \left(\begin{array}{c}x_{1}\\x_{2}\end{array}\right)
\end{array}.
$$

We start by solving the problem (\ref{s}) using the simplex method, we obtain an optimal solution equal to $x^{0}=(80,40)$. Then, we solve the following problem:
\begin{equation}\label{eqq3}
\arg \hspace{0.1cm} \underset{x\in S}\max \hspace{0.1cm} f(x)= \arg \hspace{0.1cm}\underset{x\in S}\max \hspace{0.1cm}g(x). \hspace{1.2cm}
\end{equation}
Consider the successive corners $x^{1}, x^{0}, x^{2}$ of the feasible region,  such that:
$$
\begin{array}{ccc}
x^{1}=\left(100,0\right)  & \text{and} &x^{2}=\left(60,50\right)\,.
\end{array}
$$
Then we have
 $$
\begin{array}{l}
x^{10}=\left(-20,40\right)=20\sqrt{5}\left(-\frac{1}{\sqrt{5}},\frac{2}{\sqrt{5}}\right),\vspace{0.2cm}\\
x^{02}=\left(-20,10\right)=10\sqrt{5}\left(-\frac{2}{\sqrt{5}},\frac{1}{\sqrt{5}}\right),
\end{array}
$$
and therefore we get
$$\theta_{1}=\arccos\left(-\frac{1}{\sqrt{5}}\right)\approx 116.5650^{o}$$ and $$\theta_{2}=\arccos\left(-\frac{2}{\sqrt{5}}\right)\approx 153.4349^{o}.$$ Finally, for all \hspace{0.15cm}$r>0$  and $\phi\in\left ]26.5650^{o},63.4349^{o}\right[$, the linear form:
$$
g(.)=r\left\langle (\cos(\phi),\sin(\phi)),.\right\rangle
$$
is a solution of the problem (\ref{eqq3}).\\
Finally, to do the sensitivity analysis, we write the coefficients of the objective function $f$ in the polar coordinate system, as follows:
$$
c_{0}=r_{f}\left(\cos(\phi_{f}),\sin(\phi_{f})\right),
$$
where $\phi_{f}=\arccos \left(\frac{2}{\sqrt{13}}\right)\approx 56.31^{o}$ and $r_{f}=\sqrt{13}$, and consider the solution $g$ which stabilizes $x^{0}$ as an optimal solution. After that, we calculate  $g\left(x^{0}\right)$, we get
$$
\begin{array}{lll}
g(x^{0})=g(80,40)&\approx&40\sqrt{5}r\left\langle (\cos(\phi),\sin(\phi)),(\cos(26.5650^{o}),\sin(26.5650^{o}))\right\rangle\vspace{0.25cm}\\
                 &\approx& 40\sqrt{5}r \cos(\phi-26.5650^{o})\vspace{0.2cm}\\
				&\approx& 40\sqrt{5}r\cos(\phi_{f}+\nu-26.5650^{o})\vspace{0.2cm}\\
				&\approx&40\sqrt{5}r\cos(56.31^{o}+\nu-26.5650^{o})\vspace{0.2cm}\\
&\approx& 40\sqrt{5}r\cos(\nu+29.745^{o}),
\end{array}
$$

where $\nu\in\left]-3.18^{o},33.8699^{o}\right[$ corresponds to the angle of rotation of  $\grad(f)$, to arrive at $\grad(g)$, in other words the angle between $\grad(f)$ and $\grad(g)$. So, if $\nu\in\left]-3.18^{o},0^{o}\right[$, the optimal value of $f$ increases, which means the optimal value of $g$ at $x^{0}$ is bigger than the optimal value of  $f$ at $x^{0}$. And, if $\nu\in\left]0^{o},33.8699^{o}\right[$, the optimal value of $f$ decreases.
\section{Summary.}\label{section7} In this work, we give a new formulation of the sensitivity analysis problem by the problem (\ref{Pbm3}), which allowed us to study the stability of an optimal solution of a linear programming problem in dimension 2 of the variable space. This is by reducing it to the search for a single interval by a simple calculation that does not require the use of a machine. The other advantage is that we determine all the objective functions which reach their optimum value at a common point. On the other hand, the classical approach  consists to do the study of stability under variation of just one coefficient, and the variation of the two coefficients leads us to two different intervals related to each other which makes the analysis difficult, see \cite{ref1}. Finally, we illustrated our approach with a numerical example.

 \bibliographystyle{acm}
\bibliography{references}

\end{document}